%% file: 2002-18.tex
\documentclass{gtart}

\input gtoutput

\lognumber{276}

\volumenumber{6}
\papernumber{18}
\volumeyear{2002}
\pagenumbers{523}{539}
\received{14 September 2002}
\accepted{19 November 2002}
\published{23 November 2002}
\proposed{Robion Kirby}
\seconded{Joan Birman, Vaughan Jones}

\usepackage{amsmath,amssymb,graphicx} 

\newtheorem{thm}{Theorem}[section]
\newtheorem{cor}[thm]{Corollary}
\newtheorem{lemma}[thm]{Lemma}
\newtheorem{prop}[thm]{Proposition}
\newtheorem{fact}[thm]{Fact}

\newcommand{\C}{{\mathbb C}}
\newcommand{\M}{{\mathcal M}}

\def\la{\longrightarrow}
\def\tn{\textnormal}

\def\bd{\partial}
\DeclareMathOperator{\End}{End}
\DeclareMathOperator{\Hom}{Hom}
\DeclareMathOperator{\id}{id}
\DeclareMathOperator{\tr}{tr}
\DeclareMathOperator{\genus}{genus}

\newcommand{\fig}[3]{
    \begin{figure}[ht!]
    \begin{center}\includegraphics{#2}\end{center}
    \caption{\label{#1}#3}
    \end{figure}
}

\newcommand{\figw}[4]{
    \begin{figure}[ht!]
    \begin{center}\includegraphics[width=#4]{#2}\end{center}
    \caption{\label{#1}#3}
    \end{figure}
}

\begin{document}

\title{Quantum $SU(2)$ faithfully detects mapping class\\groups modulo center}

\author{Michael H Freedman\\Kevin Walker\\Zhenghan Wang}
\shortauthors{Freedman, Walker and Wang} 
\asciiauthors{Michael H Freedman, Kevin Walker and Zhenghan Wang}

\asciiaddress{Microsoft Research, Redmond, WA 98052, 
USA (MHF and KW)\\and\\Indiana
University, Department of Mathematics\\Bloomington, IN 47405, USA}

\asciiemail{michaelf@microsoft.com, kwalker@xmission.com, zhewang@indiana.edu} 

\addresses{{\rm MHF and KW:}\qua Microsoft Research, Redmond, WA 98052, 
USA\\\smallskip\\{\rm ZW:}\qua Indiana
University, Department of Math, Bloomington, IN 47405, USA
\\\smallskip
{\tt michaelf@microsoft.com, kevin@messagetothefish.net, zhewang@indiana.edu}}

\begin{abstract}{\small
The Jones--Witten theory gives rise to representations of the
(extended) mapping class group of any closed surface $Y$ indexed by a
semi-simple Lie group $G$ and a level $k$.  In the case $G=SU(2)$
these representations (denoted $V_A(Y)$) have a particularly simple
description in terms of the Kauffman skein modules with parameter $A$
a primitive $4r^{\tn{th}}$ root of unity ($r=k+2$).  In each of these
representations (as well as the general $G$ case), Dehn twists act as
transformations of finite order, so none represents the mapping class
group $\M(Y)$ faithfully.  However, taken together, the quantum $SU(2)$
representations are faithful on non-central elements of $\M(Y)$.
(Note that $\M(Y)$ has non-trivial center only if $Y$ is a sphere with
$0, 1,$ or $2$ punctures, a torus with $0, 1,$ or $2$ punctures, or
the closed surface of genus $=2$.)  Specifically, for a non-central $h
\in \M(Y)$ there is an $r_0(h)$ such that if $r \geq r_0(h)$ and $A$
is a primitive $4r^{\tn{th}}$ root of unity then $h$ acts projectively
nontrivially on $V_A(Y)$.  Jones' \cite{J} original representation
$\rho_n$ of the braid groups $B_n$, sometimes called the generic
$q$--analog--$SU(2)$--representation, is not known to be faithful.
However, we show that any braid $h \neq \id \in B_n$ admits a cabling
$c = c_1, \ldots , c_n$ so that $\rho_N ( c(h)) \neq \id$, $N=c_1 +
\ldots + c_n$.}
\end{abstract}

\asciiabstract{The Jones-Witten theory gives rise to representations of the
(extended) mapping class group of any closed surface Y indexed by a
semi-simple Lie group G and a level k. In the case G=SU(2) these
representations (denoted V_A(Y)) have a particularly simple
description in terms of the Kauffman skein modules with parameter A a
primitive 4r-th root of unity (r=k+2). In each of these
representations (as well as the general G case), Dehn twists act as
transformations of finite order, so none represents the mapping class
group M(Y) faithfully. However, taken together, the quantum SU(2)
representations are faithful on non-central elements of M(Y). (Note
that M(Y) has non-trivial center only if Y is a sphere with 0, 1, or 2
punctures, a torus with 0, 1, or 2 punctures, or the closed surface of
genus = 2.) Specifically, for a non-central h in M(Y) there is an
r_0(h) such that if r>= r_0(h) and A is a primitive 4r-th root of
unity then h acts projectively nontrivially on V_A(Y). Jones' [J]
original representation rho_n of the braid groups B_n, sometimes
called the generic q-analog-SU(2)-representation, is not known to be
faithful.  However, we show that any braid h not= id in B_n admits a
cabling c = c_1,...,c_n so that rho_N (c(h)) not= id, N=c_1 +
... + c_n.}

\keywords{Quantum invariants, Jones--Witten theory, mapping class groups}
\asciikeywords{Quantum invariants, Jones-Witten theory, mapping class groups}

\primaryclass{57R56, 57M27}\secondaryclass{14N35, 22E46, 53D45}

\maketitlepage

\section{Introduction} \label{introsection}

Let $Y$ denote a compact, connected, oriented surface.
The
mapping class group $\M(Y) \cong \tn{Diff}^+ (Y)/\tn{Diff}^+_0
(Y)$ is defined as the orientation preserving
diffeomorphisms modulo isotopy.
(We do not put base points on
boundary components.)
Lickorish \cite{L} showed that $\M$ is finitely
generated, Hatcher and Thurston \cite{HT} showed that $\M$ is finitely
presented and explicit presentation have been written down \cite{Wj}.
It
is known that $\M$ is always residually finite \cite{G}.
Bigelow \cite{B1}
\cite{B2} has shown that $\M$ is a matrix group when $\genus(Y)=0$ and
when $Y$ is closed and $\genus(Y)=2$.
Of course, $\M(T^2 )\cong
\tn{SL}(2,Z)$ is also a matrix group.

In this note we study the quantum $SU(2)$ representations of $\M$.
Except when $\M(Y)$ is the trivial group ($Y=$ sphere or disk),
all these representations, and in fact all quantum
representations of which the authors are aware\footnote{
Bigelow's
representation is equivalent to the BMW representation but at a
generic value. At a generic value Dehn twist has infinite order
but unfortunately, generic values lead to infinite dimensional --
not quantized -- representations except in the genus $=0$ case. (To
see the difference consider admissible labelling of trees and
graphs. Even if the label set is infinite, if the labels on
valence $=1$ vertices are fixed then there are only finitely many
admissible labellings in the tree case.)
}, have kernel because Dehn
twists are carried to operators of finite order.
We prove,
however, that the direct sum of all the quantum $SU(2)$
representations is faithful except on central elements of $\M(Y)$
which are never detected.
It is well-known \cite{Iv} that $Z(\M(Y))
=\{e\}$ unless $Y= S^1 \times I, T^2 ,
T^2 -\textrm{pt}, T^2 -\textrm {2 pts}$, $T^2 \#
T^2$ in which case the center is the group generated by the
elliptic or hyper-elliptic involution.

These quantum $SU(2)$ representations are an outgrowth of
Jones--Witten theory.
We use the \cite{BHMV} construction of these
representations based on the skein theory of the Kauffman bracket.
This construction produces a projective representation $V_A(Y)$ of
$\M(Y)$ whenever Kauffman's variable $A$ is a
primitive $4r^{\tn{th}}$ root of unity.
(When $A$ is a
primitive $2r^{\tn{th}}$ root of unity a quantum--$SO(3)$
representation is the result.
All our faithfulness
results are true for this family as well.
Experts will have no
difficulty guessing the proof of this extension: simply restrict
the present proof to ``even labels''.)

First we consider surfaces $Y$ without boundary.

\begin{thm} \label{thmclosed}
Let $Y$ be a closed connected oriented surface and $\M(Y)$ its
mapping class group. For every non-central $h \in \M$, there is an
integer $r_0(h)$ such that for any $r \geq r_0 (h)$ and any $A$ a
primitive $4r^{\tn{th}}$ root of unity, the operator $\langle h
\rangle \co  V_A (Y) \to V_A(Y)$ is not the identity, $\langle h
\rangle \neq 1 \in {\mathcal P} \,\tn{End} (V_A )$, the projective
endomorphisms. In particular, any infinite direct sum of quantum
$SU(2)$ representations faithfully represents these mapping class
groups modulo center.
\end{thm}

Theorem \ref{thmclosed} and Theorem \ref{thmbdy}, which treats
surfaces with boundary, have
a formal corollary outside quantum topology (which was previously
known \cite{G}.)

\begin{cor} \label{corresid}
For all compact orientable surfaces $Y$ $\M(Y)$ is residually
finite. \qed
\end{cor}

\noindent
\proof Exploit the fact that finitely generated matrix groups
over $\C$ are residually finite.
\qed

\medskip

Within quantum topology the theorem also has an immediate
corollary.

\begin{cor} \label{cor2}
Let $Y$ be a closed connected compact orientable surface. Let $N$
be the mapping torus of a non-central $h\co Y\la Y$. Let $\langle \,
\rangle_A$ denote the closed $3$--manifold invariant associated to
$\left( SU (2), A\right), A$ a primitive $4r^{\tn{th}}$ root of
unity. For all $r \geq$ some $r_0 (h), |\langle N \rangle_A| <
|\langle S^1 \times Y\rangle_A|$.
\end{cor}

\noindent
\proof In the case of $Y$--bundles over a circle $S^1$ the
gluing relations for a TQFT imply that $\langle \, \rangle_A$ is simply
trace (monodromy) $= \tn{tr} \langle h \rangle_A$. If $\langle h \rangle_A \neq \id$ then
$|\tn{tr} \langle h \rangle_A| < |\tr \id_{V_A}|$.
\qed

\medskip

The proof of Theorem \ref{thmclosed} is relatively simple.
If $h$ is a non-central element of $\M(Y)$, then there is an embedded
curve $\alpha$ in $Y$ such that $\alpha$ and $h(\alpha)$ are not isotopic.
Associated to any curve $\alpha$ on $Y$ there is a operator
$T_\alpha\co  V_A(Y) \to V_A(Y)$, and
$T_{h(\alpha)} = \langle h \rangle T_\alpha \langle h^{-1} \rangle$.
We show that for $r$ sufficiently large $T_\alpha$ is not equal (even projectively)
to $T_{h(\alpha)}$.
It follows that $\langle h \rangle$ acts projectively nontrivially
on $V_A(Y)$.

\medskip

The rest of the paper is organized as follows.
Section \ref{qtsection} reviews the facts about the $SU(2)$ quantum invariants
we will need.
Section \ref{mtsection} contains the proofs of the main theorems, modulo
a topological lemma which is proved in Section \ref{topsection}.
Section \ref{frsection} contains further remarks on the original Jones
braid group representation.

\rk{Achknowledgements}We would like to thank Jorgen Andersen for
bringing to our attention the question of the eventual faithfulness of
the $SU(2)$ representations and for explaining to us his
gauge-theoretic approach to the problem, which he has now brought to
completion \cite{A}.  (For readers of both papers, we should point out
that it is not yet proven that the gauge theory and Kauffman bracket
constructions yield the same representations.)  We also thank the
referee for helpful comments.

Research by Wang is partially supported by NSF Grant CISE/EIA-0130388
and US Army Research Office Grant DAAD19-00-R0007.

\section{Review of $SU(2)$ quantum invariants} \label{qtsection}

In this section we briefly review Kauffman skein modules \cite{KL} and
the \cite{BHMV} construction of the $SU(2)$
quantum invariants.
For more details, see \cite{KL} and \cite{BHMV}.

The Kauffman skein module of a 3--manifold $M$ is defined to the free
vector space generated by isotopy classes of unoriented framed links in $M$,
modulo the Kauffman skein relation and replacing trivial loops with a
factor of $d = -A^2 - A^{-2}$. (See Figure \ref{fig321}.
Throughout this paper figures follow the ``blackboard framing'' convention.)
\fig{fig321}{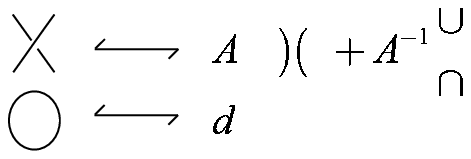}{Definition of Kauffman skein module}

One can similarly define the Kauffman skein module for a 3--manifold with
a finite collection of framed points in its boundary in terms of properly
embedded framed 1--submanifolds whose boundary is the given collection of points.
Note that for $M = S^3$ any link is equivalent to some multiple of the empty
link, so we get a $\C[A, A^{-1}]$ valued invariant of framed links on $S^3$.

In what follows we specialize to the case
\[
    A = e^{2 \pi i/4r}.
\]
(So the Kauffman ``polynomial'' of a link will actually be a complex number.)

\begin{fact} \label{fproj}
For each $k \le r-2$ there is a unique skein (finite linear combination of diagrams)
$P_k$ in $(B^3, 2k\; \mbox{points})$ such that $P_k P_k = P_k$ and $P_k$ is killed
by ``turn backs''. (See Figure \ref{fig1}.)
\end{fact}

\figw{fig1}{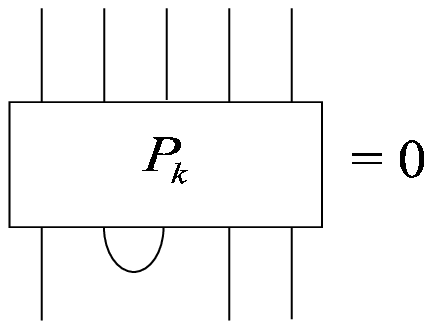}{Projector killed by turn-back}{1.3in}

It follows that $P_k$ is invariant under a 180 degree
rotation (Figure \ref{fig2}), and that $P_k$ is equal to the identity
tangle plus terms with turnbacks (Figure \ref{fig3}).
$P_k$ is called the {\it projector} on $k$ strands.

\fig{fig2}{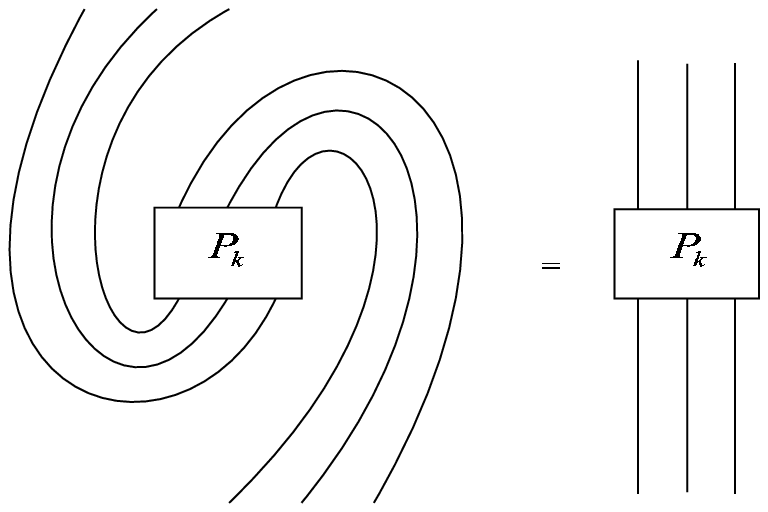}{Projector invariant under rotation}
\figw{fig3}{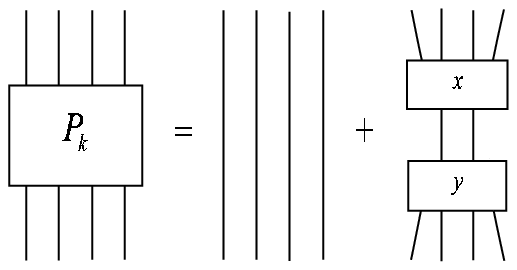}{Projector equal to identity plus turn-back terms}{2.4in}

\begin{fact} \label{frwident}
For any $n \ge 0$, then identity tangle on $n$ strands can be factored though
the sum of projectors $P_0, \ldots, P_{r-2}$ .
If $n \le r-2$, then the coefficient of $P_n$ is $1$ (Figure \ref{fig4}).
\end{fact}

\figw{fig4}{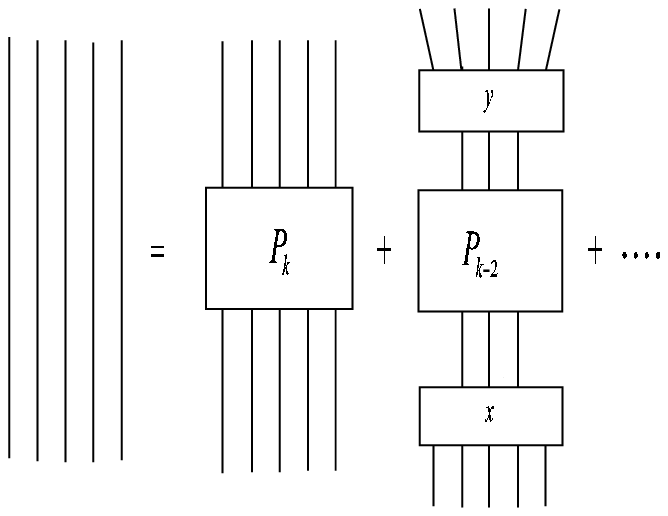}{Identity in terms of projectors}{2.7in}

The fact than only projectors up to $r-2$ are needed is a consequence of $A$ being
a $4r^{\mathrm{th}}$ root of $1$.

\begin{fact} \label{fpbraid}
Let $b$ be a braid on k strands and $c(b)$ be the signed number of crossings of $b$.
Then $b P_k = A^{c(b)} P_k$.
\end{fact}

Fact \ref{fpbraid} says that up to scalars, we can absorb a braid into a projector.
The proof follows easily from the Kauffman skein relation and Fact \ref{fproj}.

\begin{fact} \label{ftrivert}
Let $a$, $b$ and $c$ be non-negative integers
and let $X$ be a ``trivalent vertex'' skein as shown in the
left hand side of Figure \ref{fig5}.
If
(a) the three triangle inequalities are satisfied ($a \le b+c$ etc.),
(b) $a+b+c$ is even, and
(c) $a+b+c \le 2r-4$,
then $X$ is proportional to the standard diagram on the right hand side of Figure \ref{fig5}.
If these conditions are not satisfied then $X = 0$.
\end{fact}

\figw{fig5}{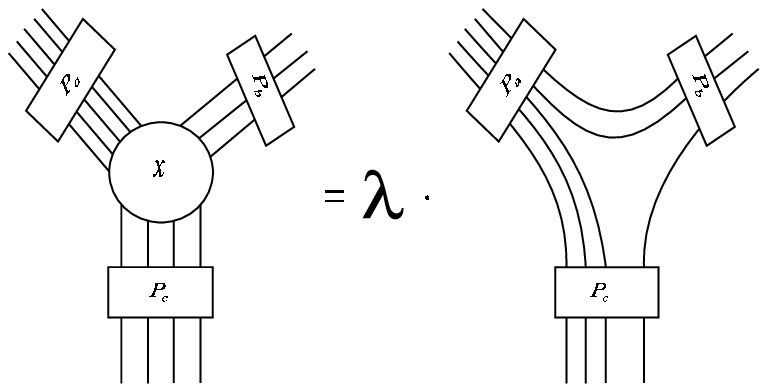}{1--dimensional trivalent vertex space}{3.2in}

Fact \ref{ftrivert} follows easily from Fact \ref{fpbraid} and Figure \ref{fig3}.


Let $G \subset M$ be a trivalent ribbon graph with edges labeled
by integers between $0$ and $r-2$, such that at each vertex the
conditions of Fact \ref{ftrivert} are satisfied. We will regard
$G$ as a shorthand notation for the linear combination of framed
links in $M$ obtained by replacing an edge of $G$ labeled by $k$
with $P_k$, and replacing trivalent vertices with the right hand
side of Figure \ref{fig5}.

Let $d_k$ be the value of the skein shown in Figure \ref{figdk} (unknot labeled
by $P_k$).
\figw{figdk}{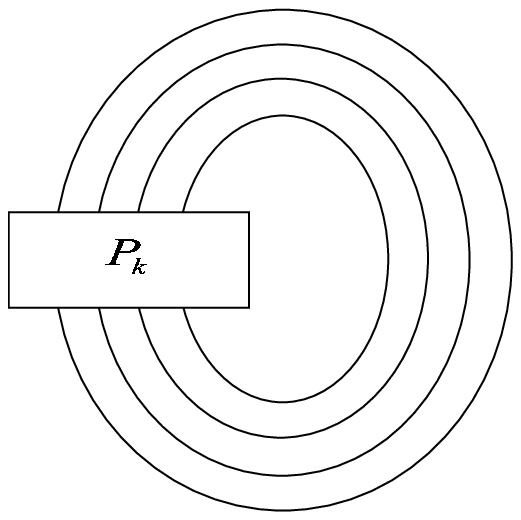}{Loop value for projector}{1.8in}
Let $s_k = c d_k$, where $c$ is a positive real number chosen so that
$\sum_{i=0}^{r-2} s_i^2 = 1$.
In a framed link diagram, a component labeled by $\omega$ will mean
the linear combination shown in Figure \ref{figomdef}.

\fig{figomdef}{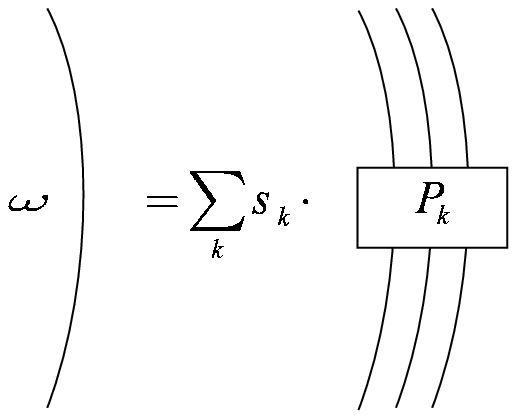}{Definition of $\omega$ label}

\begin{fact} \label{fsurg}
Framed links with components labeled by $\omega$ are invariant under handle slides,
balanced stabilization, and the introduction of a circumcision pair.
(See Figures \ref{fighs}, \ref{figbstab} and \ref{figcpair}.)
\end{fact}

\figw{fighs}{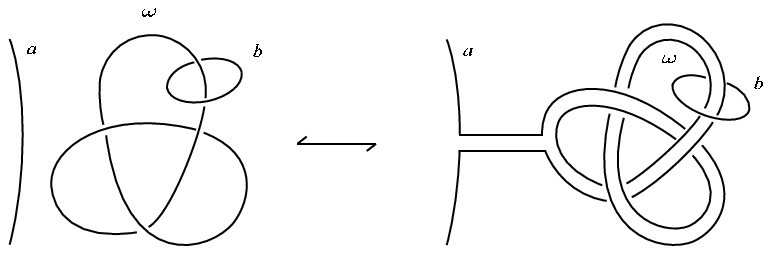}{Handle slide invariance}{4in}
\fig{figbstab}{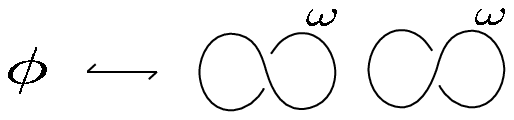}{Balanced stabilization invariance}
\figw{figcpair}{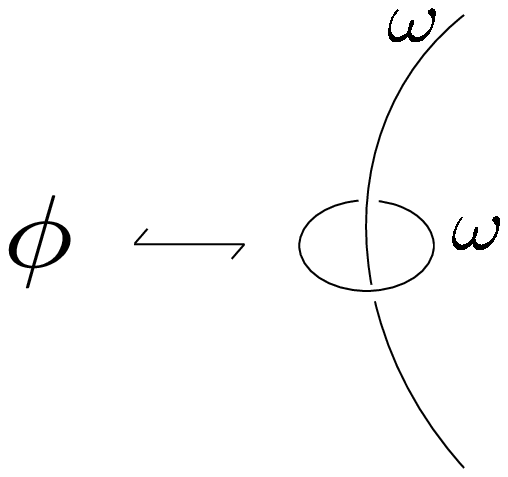}{Circumcision pair invariance}{1.3in}

Let $L$ be a framed link in $S^3$. Let $L_\omega$ be the linear
combination of labeled framed links obtained by labelling each
component of $L$ by $\omega$. It follows from Fact \ref{fsurg}
that the Kauffman polynomial of $L_\omega$ depends only on the
3--manifold described by interpreting $L$ as a surgery diagram, and
on the signature of $L$. For any closed, oriented 3--manifold $M$
and integer $n$ define $Z(M,n)$ to be this invariant (ie, $Z(M,
n)$ is equal to the Kauffman polynomial of $L_\omega$, where $L
\subset S^3$ is any surgery description of $M$ with signature
$n$.) It is easy to see that $Z(M, n) = C^{n-m} Z(M, m)$, where
$C$ is the value of the Kauffman polynomial of an unknot with
framing $1$ (right handed twist).

(Note: $n$ can be interpreted as an equivalence class of framings of the tangent
bundle of $M$, a bordism class of null-bordisms of $M$, or a $p_1$--structure on $M$.
See \cite{A}, \cite{Wa} and \cite{BHMV}.)

Next we follow the \cite{BHMV} approach to construct a vector space $V(Y)$ for each
closed, oriented 2--manifold $Y$, and an
invariant $Z(M) \in V(\bd M)$ for an oriented 3--manifold
with boundary.
These 2--manifolds and 3--manifolds with boundary should also be equipped
with extra structure (framing, null-bordism, or $p_1$--structure), but
we will suppress mention of this since the arguments in the remainder
of the paper work even with a projective ambiguity.

Let $Y$ be a closed, oriented 2--manifold.
Let $\bd^{-1}Y$ be the set of all isomorphism classes of
pairs $(M, L)$, where $\bd M = Y$ and $L$ is a labeled ribbon graph in the interior
of $M$.
Let W(Y) be the free vector space generated by $\bd^{-1}Y$.
There is a pairing $W(Y) \otimes W(-Y) \to \C$ given by
$x \otimes y \mapsto Z(x \cup y)$.
Define $V(Y)$ to be the quotient of $W(Y)$ by the annihilator of $W(-Y)$ with respect
to this pairing.
In other words, $x \sim x'$ if $Z(x \cup y) = Z(x' \cup y)$ for all $y \in W(-Y)$.

If $Y$ is not closed choose a labelling $l$ of the boundary components of $Y$
by integers $0 \le l_c \le r-2$.
Let $\widehat{Y}$ be the result of capping off each boundary component of $Y$
by $D^2$.
Define $\partial^{-1}(Y, l)$ to be the set of isomorphism classes of 3--manifold $M$
with $\partial M$ identified with $\widehat{Y}$,
and with a properly embedded framed tangle in $M$ which coincides with a standardly
embedded copy of $P_k$ in a collar neighborhood of each cap disk, where $k$ is the
label assigned to that boundary component of $Y$ by $l$.
We can now define $V(Y; l)$ as above.

The extended mapping class group of $Y$ acts on $\bd^{-1}Y$, and thus on
$V(Y)$.
The ordinary, non-extended mapping class group of $Y$ has a projective action
on $V(Y)$.

The surgery formula for $Z$ shows that $V(Y)$ is spanned by the equivalence classes
of links in any single 3--manifold $M$, $\bd M = Y$.
For example, we could take $M$ to be a handlebody $H$ (assuming $Y$ is connected).
It then follows from Facts \ref{frwident} and \ref{ftrivert} that:

\begin{fact} \label{fhbbasis}
Let $H$ be a handlebody with spine $S$,
(ie, $S$ is a 1--complex with vertices at most trivalent, and $H$ is a regular
neighborhood of $S$.)
Then $V(\bd H)$ has a basis corresponding to all labellings of the 1--cells
of $S$ by integers between $0$ and $r-2$, such that the parity and quantum triangle inequalities
of Fact \ref{ftrivert} are satisfied at each vertex of $S$.
\end{fact}

If $Y$ has non-empty boundary, we get a basis of $V(Y,l)$
by letting $\widehat{Y}$ bound a handlebody $H$ and considering spines of $H$
which meet each cap disk of $\widehat{Y}$ once.
Labellings of the spine are constrained to agree with $l$ on 1--cells meeting the boundary.

If $Y$ is closed then $\End(V(Y))$ can be identified with $V(Y \coprod -Y)$, and so is spanned
by elements of the form $Z(Y \times I, L)$, where $L$ is a labeled framed link in $Y \times I$.
If $Y$ has boundary then $\bigoplus_l \End(V(Y, l))$ can be identified with $V(D(Y))$,
where $l$ runs through all labellings of $\bd Y$ and $D(Y) = Y \cup_{\bd Y} -Y$ is the double
of $Y$ along its boundary.
$D(Y)$ bounds $Y \times I$, and as before $\bigoplus_l \End(V(Y, l))$ is
spanned by elements of the form $Z(Y \times I, L)$, where $L$ is a
labeled framed link in $Y \times I$.
In both cases the action of $\End(\ldots)$ is given in geometric terms by gluing
$(Y \times I, L)$ onto a 3--manifold (bounded by $Y$) representing
an element of $V(Y)$ (or $V(Y,l)$).

\section{Proof of main theorems} \label{mtsection}

Let $Y$ be a closed, oriented surface, $h\co Y \to Y$ an orientation preserving
homeomorphism, and $V_h\co V(Y) \to V(Y)$ the action of $h$ on the TQFT vector space.

\begin{prop} \label{mpclosed}
Suppose there exists an unoriented simple closed curve $a \subset Y$ such that
$h(a)$ is not isotopic (as a set) to $a$.
Then $V_h$ is a multiple of the identity for at most finitely many $r$.
That is, as $r$ increases $h$ is eventually detected.
\end{prop}

\noindent
\proof
Let $C(a) = Z(Y \times I, a \times \{1/2\}) \in V(Y) \otimes V(-Y) = \End(V(Y))$.
Define $C(h(a))$ similarly.
It's easy to see that $C(h(a)) = V_h C(a) V_h^{-1}$.
It therefore suffices to show that $C(a) \neq C(h(a))$.

By Lemma \ref{closedtl} there exists a handlebody $H$ bounded by $Y$ such that
$a$ bounds an embedded disk in $H$ and $h(a)$ is a non-trivial ``graph geodesic''
with respect to a spine $S$ of $H$.
Let $Z(H) \in V(Y)$ be the vector determined by $H$, and
$Z(H, h(a)) \in V(Y)$ be the vector determined by the pair $(H, h(a))$.
(We can push $h(a)$ into the interior or $H$.)
Then
\begin{gather*}
    C(a)(Z(H)) = Z(H, a) = d \cdot Z(H) ,\\
{\rm and}\qquad\qquad    C(h(a))(Z(H)) = Z(H, h(a)) .
\qquad\qquad\phantom{\rm and}\end{gather*}
It therefore suffices to show that $Z(H, h(a))$ is not a multiple of $Z(H)$.

For each edge $e$ of the spine $S$, let $w_e$ be the (unsigned) number of times $h(a)$
passes over $e$.
Let $m$ be the maximum of all $w_e+w_f+w_g$ such that $e$, $f$ and $g$ meet
at a vertex of $S$.
Choose $r$ such that $2r-4 \ge m$.

Let $b_w$ be the basis vector of $V(Y)$ corresponding the labelling $w$.
We claim that $Z(H, h(a)) = \lambda b_w + v$, where $\lambda \ne 0$ and $v$ consists of
``lower order'' terms -- multiples of $b_v$, where $v_e \le w_e$ for all edges $e$ of $S$
and $v \ne w$.
This follows from Facts \ref{frwident}, \ref{ftrivert} and \ref{fpbraid}.
Apply Fact \ref{frwident} at each edge of $S$.
Apply Fact \ref{ftrivert} at each vertex to see that the result is a linear
combination of $b_w$ and lower order terms.
Fact \ref{fpbraid} and the graph geodesic property of $h(a)$ show that the coefficient
of $b_w$ is non-zero.
On the other hand, $Z(H)$ is the basis vector corresponding to the zero (empty)
labelling of $S$.
\qed

\medskip

\noindent
\proof[Proof of Theorem \ref{thmclosed}]
By Lemma \ref{mcgcenter}, non-central elements of the mapping
class group must move a simple closed curve, so Theorem \ref{thmclosed}
follows from Proposition \ref{mpclosed}.
\qed

\medskip

Next we consider the case where $Y$ has boundary.
As before, let $h\co Y \to Y$ be an orientation preserving
homeomorphism and
\[
    V_h \in \bigoplus_{l,l'} \Hom(V(Y,l), V(Y,l'))
\]
be the action of $h$ on the TQFT vector spaces.

\begin{prop} \label{mpbdy}
Suppose there exists an unoriented, homologically essential
simple closed curve $a \subset Y$ such that
$h(a)$ is not isotopic to $a$.
Then $V_h$ is a multiple of the identity for at most finitely many $r$.
That is, as $r$ increases $h$ is eventually detected.
\end{prop}

\noindent
\proof
Define operators $C(a)$ and $C(h(a))$ as in the proof of Proposition \ref{mpclosed}.
(Note that while $V_h \in \bigoplus_{l,l'} \Hom(V(Y,l), V(Y,l'))$,
$C(a)$ and $C(h(a))$ lie in the block diagonal $\bigoplus_l \End(V(Y,l)$.)
As before, it suffices to show that $C(a) \ne C(h(a))$.

By Lemma \ref{bdyspine}, $a \times \{1/2\}$ can be extended to a
spine of $Y \times I$. Since $h(a)$ is not isotopic in $Y$ to $a$,
$h(a)$ must be isotopic to a graph geodesic distinct from $a
\times \{1/2\}$. It follows from Fact \ref{fhbbasis} that $C(a)$
and $C(h(a))$ are (projectively) distinct elements in $V(\bd(Y
\times I)) = \bigoplus_l \End(V(Y,l))$, provided $r$ is
sufficiently large. \qed

\medskip

We can now prove:

\begin{thm} \label{thmbdy}
Let $Y$ be a connected orientable surface with boundary and let
$h$ be a non-central diffeomorphism of $Y$. Let $V_h \in
\bigoplus_{l,l'} \Hom(V(Y,l), V(Y,l'))$ be the action of $h$ on
the TQFT vector spaces. Then $V_h$ is a multiple of the identity
for at most finitely many $r$.
\end{thm}

\proof In light of Proposition \ref{mpbdy}, it
suffices to show that any diffeomorphism of $Y$ which fixes all
homologically essential simple closed curves lies in the center of
the mapping class group. Let $h$ be such a diffeomorphism. Then
unless $Y$ is
an
annulus $h$ cannot permute the boundary components
of $Y$; also $h$ commutes with Dehn twists along homologically
essential curves and all $\lq\lq$essential" braid twists $b$
($1/2$ Dehn twists
which permute a pair of boundary components)
along an essential scc $\gamma$ which bounds
a pair of pants to at least one side. Letting $\M(Y)$ denote the
full mapping class group and $N$ the number of boundary components
of $Y$ we have a short exact sequence:
\[
1 \rightarrow \M_0 (Y) \rightarrow \M(Y) \rightarrow \sigma
(N) \rightarrow 1
\]
where $\sigma(N)$ is the permutation group and $\M_0 (Y)$ the
kernel. If $N=1$, $\M(Y)= \M_0(Y)$ is generated by Dehn
twists along essential sccs and if $N\geq 3$, $\M(Y)$ is generated
by Dehn twists along essential sccs together with essential braid
twists $b$ as above. In these cases $h$ commutes with a
generating set, and therefore all, of $\M(Y)$. When $N=2$ we
need to include some (any) $\lq\lq$inessential" braid twist $b'$
along a scc $\gamma'$ bounding a pair of pants on one side and null
bounding on the other side. Since $\gamma'$ is null homologous,
special pleading is now required to prove that $h(\gamma')\approx
\gamma'$. We exploit the fact that we may pick any $\gamma'$ we
like so long as it cobounds a pair of points with $\partial Y$.
Choosing $\gamma'$ amounts to picking a simple arc $\alpha$
between the two components $\partial^{+}$ and $\partial^{-}$ of
$\partial Y$ (and then thickening). Choose $\alpha$ so that the
geometric intersection numbers are $(\alpha , \beta_0 )=1,
(\alpha, \beta_1 )=0, \cdots , (\alpha , \beta_{2g} ) =0$, where
$\{\beta_0 , \ldots , \beta_{2g}\}$ is a chain a $2$ genus
$(Y)+1$ sccs in int$(Y)$ so that only $\beta$'s of adjacent
indices meet and these meet transversely in a single point and so
that $\partial^{+}$ is separated from $\partial^{-}$ by
$\overset{2g}{\underset{i=0}{U}}\,\beta_i$ (see Figure 11). Now
$\left(h(\alpha) , \beta_i\right) = \left(h(\alpha),
h(\beta_i)\right) = (\alpha , \beta_i) = \delta(i)$. It follows
that $h(\alpha)$ is isotopic back to $\alpha$ (The isotopy may
twist $\partial Y$.) and that $h(\gamma')\approx\gamma'$. Now the
proof can be finished for $N=2$, as in the case $N\geq3$, by
taking a generating set for $\M(Y)$ consisting of $b' = b'
(\gamma')$ together with Dehn twists about essential sccs. \qed

{\nocolon\figw{fig_11}{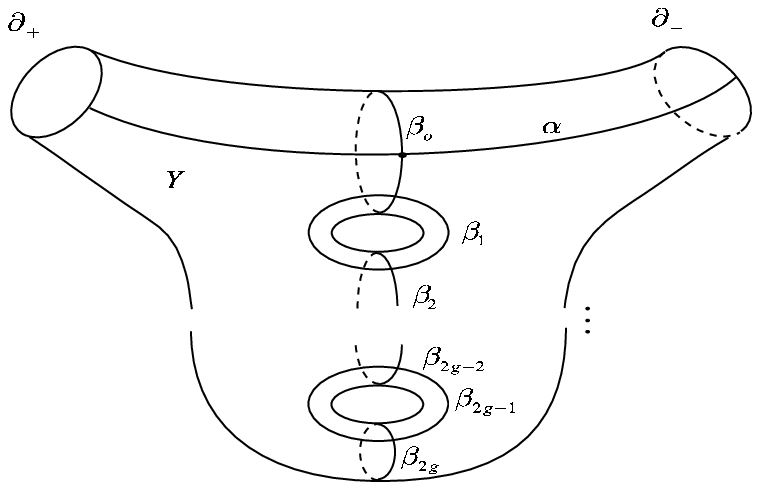}{}{3.5in}}

\section{Some topological lemmas} \label{topsection}

For applications to closed surfaces, we need:

\begin{lemma} \label{closedtl}
Let $a$ and $b$ be two non-trivial, non-isotopic simple closed
curves on a closed orientable surface $Y$. Then there exists a
pants decomposition of $Y$ such that $a$ is one of the decomposing
curves and $b$ is a non-trivial ``graph geodesic'' with respect to
the decomposition. (That is, $b$ does not intersect any curve of
the decomposition twice in a row.)
\end{lemma}

\proof
We will inductively choose a set of decomposing curves on $Y$, starting with $a$.
At each stage, let $Y'$ denote $Y$ cut along the curves we have chosen thus far, and let $b'$
denote the image of $b$ in $Y'$.
$b'$ is a properly embedded, possibly disconnected, 1--submanifold of $Y'$.

We say that $Y'$ and $b'$ satisfy Condition X if for each component $S$ of
$Y'$ and each component $e$ of $S \cap b'$ either (a) $e$ is non-separating or
(b) each component of $S \setminus e$ has genus greater than zero.

Note that initially, when $Y'$ is $Y \setminus a$, Condition X is satisfied (after possibly
isotoping $b$ to remove bigons with $a$).
If $Y'$ consists only of pairs of pants
(or an annulus if $Y$ was a torus), then Condition X implies the
graph geodesic property.
Thus it suffices to show that at each stage we can choose an additional
decomposing curve such that
Condition X is preserved, until we have a pants decomposition.

Choose a component $S$ of $Y'$ which is not a pair of pants or annulus.
We will find a simple closed curve (scc) $c$ in $S$ such that $S \setminus c$
still satisfies condition X.

If $S$ has genus greater than zero, let $\bar{S}$ be the closed
surface obtained by capping of the boundary of $S$ with disks.
Those components of $b' \cap S$ which are: (1) an arc with both
endpoints on the same boundary component of $S$, or (2) a scc,
determine a well-defined isotopy class of curves in $\bar{S}$. In
case (1) complete the arc to a circle by coning its endpoints in
the cap; in case (2) simply include. Choose a curve $c$ in $S$
whose image in $\bar{S}$ does not lie in any of the aforementioned
isotopy classes.    
If the genus of $S$ is $\ge 2$, we further require that $c$ is a 
separating curve.
By pushing $c$ across punctured bigons, we may
assume that no component of $S \setminus (c \cup b')$ is a
punctured bigon (see Figure \ref{figbbigon}). Thus Condition X is
satisfied.

\figw{figbbigon}{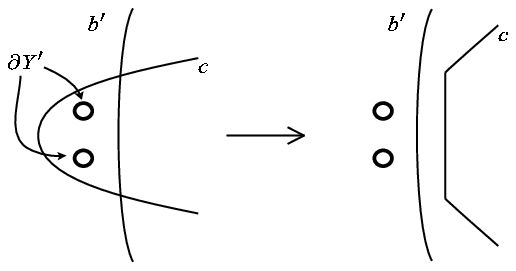}{Push across punctured bigon}{3in}

Note that for a genus 0 surface, Condition X is satisfied if and
only if all components of $b'$ are arcs which connect distinct
boundary components. Assuming $S$ has four or more punctures, we
need to find a scc $c \subset S'$ which is not boundary parallel
and meets each arc of $b'$ in at most one point. Cutting along
$c$ perpetuates condition $X$. We use a little geometry here to
avoid a greater amount of combinatorics. A well known theorem of
K$\mathrm{\ddot{o}}$ebe\footnote{Often called Andreev's Theorem.}\cite{K}
 represents the edges of any spherical graph by disjoint geodesic
arcs of length $<\pi$. Regarding the punctures of $S$ as vertices,
represent $b'$ in this way, with the understanding that parallel
arcs of $b'$ collapse to a single edge. We call two arcs of $b'$
parallel if they join the same boundary components $x$ and $y$,
and together with an arc in $x$ and an arc in $y$, bound a
rectangle in $S$. Any great circle $\gamma$ disjoint from the
vertices and containing at least two vertices in each
complementary hemisphere is a good choice for $c$. To find such a
$\gamma$, start with the great circle $\gamma'$ determined by any
two nonantipotal vertices and perturb it suitably.\qed


\begin{lemma} \label{bdyspine}
Let $Y$ be a connected orientable surface with boundary and let
$a$ be a homologically essential simple closed curve in $Y$. Then
$a$ can be extended to a spine of $Y$.
\end{lemma}

\proof
Cut $Y$ along $a$ and use the classification of surfaces.
\qed

\begin{lemma} \label{mcgcenter}
Suppose $Y$ is a compact oriented surface
with or without boundary.
Suppose $h\co  Y \la Y$ is an orientation
preserving homeomorphism, not isotopic to $id_{Y}$, which does not
change the unparameterised isotopy class of any scc in $Y$.
Then $Y$ is either an annulus, a torus, a torus with $\leq 2$ punctures, or the
closed surface of genus $=2$ and $h$ is either the elliptic or
hyperelliptic involution.
\end{lemma}

\proof
If $h\co  Y \la Y$ leaves all (unoriented) isotopy
classes of scc's invariant then $h$ will commute with all Dehn
twists.
Since Dehn twists generate $\M(Y)$ \cite{L}, $h \in
(\tn{Center} (\M(Y)) =: Z (\M(Y))$.
It is well-known (\cite{Iv}, Theorem
7.5D) that the only surfaces with $Z (\M(Y)) \neq \{ e\}$ are $Y =
T^2$, $T^2 \setminus$ pt., $T^2 \setminus$ 2 pts. , $S^1 \times I$, and $T^2 \#
T^2$. Furthermore the only nontrivial element of these centers are
the elliptic and hyperelliptic involutions respectively.
\qed

\section{Further remarks} \label{frsection}

The Jones representation contains the Burau representation as a
particular summand. It is known that the Burau representation is
not faithful for $B_n$ with $n\geq 5$. On the other hand, the
Jones representation can be obtained by specializing the BMW
representation which is faithful for $B_n$. It seems hard to
decide the faithfulness of Jones representation but in this
direction, we prove:

\begin{thm} \label{jonesthm}
For every braid $h
\neq 1 \in B_n$, the $n$--strand (unframed) braid group $n\geq 2$, there is a
cabling $(c_1 , \ldots , c_n)$ of $h$ on which the $SU(2)$--Jones
representation is nontrivial.
\end{thm}

\proof
The Jones representation on $B_n$ when specialized to
$A = e^{2 \pi i/ 4r}$, $t = e^{2 \pi i /r}$ decomposes as a
direct sum of singular and nonsingular pieces.
The nonsingular
piece is a  sum of the
$SU(2)$--quantum representations on $V_{1_1 \ldots , 1_n , m}$,
the Hilbert space at level $k = r-2$.
The subscripts of $V$ are
admissible labels at finite punctures and infinity.
Cabling
produces sums of irreducibles according to a Clebsch--Gordon
formula.
In particular the Jones representation on the $c_1 ,
\ldots , c_n$ cabling contains as a summand a copy of each
admissible $V_{c_1 , \ldots , c_n , m}$.
Thus it is sufficient to
prove that $h$ acts nontrivially on at least one of these.
Theorem \ref{thmbdy} says that with only finitely many exceptions $h$ is nontrivial
in these representations, provided $h$ is not homotopic to the
identity in the $n+1 -$ punctured sphere, that is $[h] \neq 1 \in
(\tn{spherical braid group})_{n+1} = SB_{n+1}$.

The proof is not yet finished since the natural morphism $B_n \la
SB_{n+1}$ has kennel = center$(B_n ) = \langle \tn{full
twist}\rangle$. This $\lq\lq$full twist" is Dehn twist about
infinity and although this twist is trivial in all End$((V_{c_1 ,
\ldots , c_n , m})$ its action is computed \cite{KL} to be
multiplication by the unit scalar $A^{m (m+2)}$. Thus each nontrivial central
element is also detected in infinitely many $V$'s. \qed

\end{document}

%% file: gtoutput.tex
\def\ifplaintex{\expandafter\ifx\csname documentclass\endcsname\relax}

\ifplaintex 
\else
\fi

\expandafter\ifx\csname beginpicture\endcsname\relax
\expandafter\ifx\csname documentclass\endcsname\relax
\input pictex \else
\input prepictex \input pictex \input postpictex \fi\fi

\def\gt{{\mathsurround=0pt\it $\cal G\mskip-2mu$eometry \&\ 
$\cal T\!\!$opology}}        

\def\gtp{{\mathsurround=0pt\it $\cal G\mskip-2mu$eometry \&\ 
$\cal T\!\!$opology $\cal P\!$ublications}}  

\def\lognumber#1{\def\thelognumber{#1}}
\def\volumenumber#1{\def\thevolumenumber{#1}}
\def\papernumber#1{\def\thepapernumber{#1}}
\def\volumeyear#1{\def\thevolumeyear{#1}}

\def\pagenumbers#1#2{\def\startpage{#1}\def\finishpage{#2}}
\def\published#1{\def\publishdate{#1}}
\def\proposed#1{\def\theproposer{#1}}
\def\seconded#1{\def\theseconders{#1}}
\def\received#1{\def\receiveddate{#1}}

\def\accepted#1{\def\accepteddate{#1}}

\def\asciiauthors#1{\def\theasciiauthors{#1}}
\def\asciiaddress#1{\def\theasciiaddress{#1}}
\def\asciiemail#1{\def\theasciiemail{#1}}
\long\def\asciiabstract#1{\long\def\theasciiabstract{#1}}
\def\asciikeywords#1{\def\theasciikeywords{#1}}
\def\shortauthors#1{\def\theshortauthors{#1}}

\let\\\par\let\thelognumber\relax
\let\thevolumenumber\relax\let\thepapernumber\relax
\let\thevolumeyear\relax\let\thesamplenumber\relax\let\startpage\relax
\let\finishpage\relax\let\publishdate\relax\let\receiveddate\relax
\let\reviseddate\relax\let\accepteddate\relax\let\theasciititle\relax
\let\theasciiauthors\relax\let\theasciiaddress\relax
\let\theasciiabstract\relax\let\theasciikeywords\relax
\let\theasciiemail\relax\let\theshortauthors\relax\let\theshorttitle\relax

\long\def\maketitlep{   

\count0=\startpage

\gt\hfill      
\beginpicture
\setcoordinatesystem units <0.33truein, 0.33truein> point at 2.2 0.9
\setplotsymbol ({$\cal G$})
\plotsymbolspacing=9truept
\circulararc 315 degrees from 0 1 center at 0 0
\setplotsymbol ({$\cal T$})
\circulararc 315 degrees from 1 -1 center at 1 0
\endpicture
\break
{\small\ifx\thesamplenumber\relax 
Volume \else Sample
\fi\thevolumenumber\ (\thevolumeyear)
\startpage--\finishpage\nl
Published: \publishdate}
\vglue 0.5truein plus 0.4fil minus 0.1truein

{\parskip=0pt\leftskip 0pt plus 1fil\def\\{\par\smallskip}{\ifplaintex\large
\else\Large\fi\bf\thetitle}\par\medskip}   

\vglue 0pt plus 0.1fil 

{\parskip=0pt\leftskip 0pt plus 1fil\def\\{\par}{\sc\theauthors}
\par\medskip}

\vglue 0pt plus 0.1fil 

{\small\parskip=0pt\let\newline\\
{\leftskip 0pt plus 1fil\def\\{\par}{\sl\theaddress}\par}
\expandafter\ifx\theemail\relax    
\relax\else\vglue 5pt plus 0.02fil minus 2pt\def\\{\stdspace{\rm 
and}\stdspace} 
\cl{Email:\stdspace\tt\theemail}\fi
\ifx\theurl\relax                  
\relax\else\vglue 5pt plus 0.02fil minus 2pt\def\\{\stdspace{\rm 
and}\stdspace}
\cl{URL:\stdspace\tt\theurl}\fi\par}

\vglue 7pt plus 0.3fil minus 3pt

{\bf Abstract}
\vglue 5pt plus 0.1fil minus 2pt

\theabstract

\vglue 7pt plus 0.3fil minus 3pt

{\bf AMS Classification numbers}\quad Primary:\quad \theprimaryclass

Secondary:\quad \thesecondaryclass

\vglue 5pt plus 0.3fil minus 2pt

{\bf Keywords}\quad \thekeywords

\vglue 10pt plus 0.5fil minus 5pt

{\small  Proposed: \theproposer\hfill Received: \receiveddate\nl
Seconded: \theseconders\hfill 
\ifx\reviseddate\relax                         
Accepted: \accepteddate                        
\else
Revised: \reviseddate                          
\fi}
\eject
}       

\let\maketitlepage\maketitlep
\let\maketitle\maketitlepage

\font\phead=cmsl9 scaled 950
\font=cmsl9 scaled 1050
\font\pnum=cmbx10 scaled 913
\font\lnum=cmbx10 
\font\pfoot=cmsl9 scaled 950
\font=cmsl9 scaled 1050
\ifplaintex
\headline{\vbox to 0pt{\vskip -4.5mm\line{\small\phead\ifnum
\count0=\startpage ISSN 1364-0380 (on line)
1465-3060 (printed) \hfill {\pnum\folio}\else\ifodd\count0\def\\{ }%
\ifx\theshorttitle\relax\thetitle\else\theshorttitle\fi\hfill{\pnum\folio}
\else\def\\{ and }{\pnum\folio}\hfill\ifx\theshortauthors\relax\theauthors
\else\theshortauthors\fi\fi\fi}\vss}}
\footline{\vbox to 0pt{\vglue 0mm\line{\small\pfoot\ifnum\count0=\startpage
\copyright\ \gtp\hfill\else
\gt, Volume \thevolumenumber\ (\thevolumeyear)\hfill\fi}\vss
}}
\else
\makeatletter
\def\@oddhead{{\small\lhead\ifnum\count0=\startpage ISSN 1364-0380 (on line)
1465-3060 (printed) \hfill {\lnum\number\count0}\else\ifodd\count0
\def\\{ }\ifx\theshorttitle\relax \thetitle \else\theshorttitle\fi\hfill
{\lnum\number\count0}\else\def\\{ and }{\lnum\number\count0}
\hfill\ifx\theshortauthors\relax 
\theauthors\else\theshortauthors\fi\fi\fi}}\def\@evenhead{\@oddhead}
\def\@oddfoot{\small\lfoot\ifnum\count0=\startpage\copyright\ \gtp\hfill\else
\gt, Volume \thevolumenumber\ (\thevolumeyear)\hfill\fi}
\def\@evenfoot{\@oddfoot}
\makeatother
\fi

\newwrite\gtoutfile
\long\gdef\makeheadfile{  
{\def\\{, }\def\s{ }
\immediate\openout\gtoutfile head.xxx
\immediate\write\gtoutfile{To: math@arxiv.org}
\immediate\write\gtoutfile{Subject: put or rep NNNNN:pppp}
\immediate\write\gtoutfile{--text follows this line--}
\immediate\write\gtoutfile{Proxy-for: \ifx\theasciiauthors\relax
\theauthors\else\theasciiauthors\fi\s<\ifx\theasciiemail\relax\theemail\else\theasciiemail\fi>}
\immediate\write\gtoutfile{\noexpand\\}
\immediate\write\gtoutfile{Authors: \ifx\theasciiauthors\relax
\theauthors\else\theasciiauthors\fi}
{\def\\{ }\immediate\write\gtoutfile{Title: \ifx\theasciititle\relax
\thetitle\else\theasciititle\fi}}
\immediate\write\gtoutfile{Subj-class: GT or SG or MG etc}
\immediate\write\gtoutfile{MSC-class: \theprimaryclass\ifx\thesecondaryclass\relax\else, \thesecondaryclass\fi}
\immediate\write\gtoutfile{Journal-ref: Geom. Topol. \thevolumenumber
(\thevolumeyear) \startpage-\finishpage}
\immediate\write\gtoutfile{Comments: Published by Geometry and Topology at}
\immediate\write\gtoutfile{\s\s http://www.maths.warwick.ac.uk/gt/GTVol\thevolumenumber/paper\thepapernumber.abs.html}
\immediate\write\gtoutfile{\noexpand\\}
\immediate\write\gtoutfile{}
\ifx\theasciiabstract\relax
\immediate\write\gtoutfile{\theabstract}\else
\immediate\write\gtoutfile{\theasciiabstract}\fi
\immediate\write\gtoutfile{}
\immediate\write\gtoutfile{\noexpand\\}
\immediate\write\gtoutfile{}
\immediate\closeout\gtoutfile}}  

\def\maketitlepage{\maketitlep\makeheadfile}
\let\maketitle\maketitlepage